\RequirePackage{etoolbox}
\csdef{input@path}{{style/}{graphics/}}
\documentclass[MSNbibl,number,citesort,seceqn,dvips]{arxbj}
\usepackage{upgreek}

\aid{0}
\volume{21}
\issue{1}
\pubyear{2015}
\firstpage{489}
\lastpage{504}
\doi{10.3150/13-BEJ576} 

\makeatletter

\renewcommand{\emptyset}{\varnothing}
\renewcommand{\pi}{\uppi}

\def\E{{\mathbb{E}}}
\def\N{{\mathbb{N}}}
\def\P{{\mathbb{P}}}
\def\R{{\mathbb{R}}}
\def\T{{\mathcal{T}}}

\newcommand{\wt}{\widetilde}

\def\N{\mathbb{N}}
\def\E{\mathbb{E}}
\def\P{\mathbb{P}}

\renewcommand{\d}{\delta}
\renewcommand{\a}{\alpha}
\renewcommand{\b}{\beta}
\newcommand{\lam}{\lambda}
\newcommand{\eps}{\varepsilon}
\newcommand{\g}{\gamma}
\newcommand{\ov}{\overline}
\newcommand{\wh}{\widehat}
\newtheorem{theorem}{Theorem}[section]
\newtheorem{cor}[theorem]{Corollary}
\newtheorem{lem}[theorem]{Lemma}
\newtheorem{prop}[theorem]{Proposition}
\newremark{rem}[theorem]{Remark}

\newcommand{\eqref}[1]{(\ref{#1})}

\def\sfrac#1#2{#1/#2}

\def\afrac#1#2{#1/(#2)}

\def\sklfrac#1#2{(#1/#2)}
\makeatother

\begin{document}
\begin{frontmatter}

\title{Precise tail asymptotics of fixed points of the smoothing
transform with general weights}
\runtitle{Tail asymptotics of fixed points of the smoothing transform}

\begin{aug}
\author[a]{\inits{D.}\fnms{D.} \snm{Buraczewski}\corref{}\thanksref{e1}\ead[label=e1,mark]{dbura@math.uni.wroc.pl}},
\author[a]{\inits{E.}\fnms{E.} \snm{Damek}\thanksref{e2}\ead[label=e2,mark]{edamek@math.uni.wroc.pl}} 
\and
\author[a]{\inits{J.}\fnms{J.} \snm{Zienkiewicz}\thanksref{e3}\ead[label=e3,mark]{zenek@math.uni.wroc.pl}}
\address[a]{Instytut Matematyczny, Uniwersytet Wroclawski, 50-384 Wroclaw,
pl. Grunwaldzki 2/4, Poland.\\ \printead{e1,e2,e3}}
\end{aug}

\received{\smonth{10} \syear{2012}}
\revised{\smonth{10} \syear{2013}}

%
\begin{abstract}
We consider solutions of the stochastic equation $R=_d
\sum_{i=1}^N A_iR_i + B$, where $N>1$ is a fixed constant, $A_i$
are independent, identically distributed random variables and
$R_i$ are independent copies of $R$, which are independent
both from $A_i$'s and $B$. The hypotheses ensuring existence of
solutions are well known. Moreover under a number of assumptions
the main being $\E|A_1|^\a= 1/N$ and $\E|A_1|^\a\log|A_1| >0$,
the limit $\lim_{t\to\infty } t^\a\P[|R|>t] = K$ exists. In the
present paper, we prove positivity of $K$.
\end{abstract}

%

\end{frontmatter}

\section{Introduction}

Let $N>1$ be an integer, $A_1,\ldots,A_N,B$ real valued
random variables such that $A_i$ are independent and identically
distributed (i.i.d.). On the set $P (\R)$ of
probability measures on the real line the smoothing transform is
defined as follows
\[
\mu\mapsto{\mathcal L} \Biggl(\sum_{i=1}^NA_iR_i
+ B \Biggr),
\]
where $R_1,\ldots,R_N$, are i.i.d. random variables with common
distribution $\mu$, independent of $(B, A_1,\ldots, A_N)$ and
${\mathcal L} (R)$ denotes the law of the random variable $R$. A
fixed point of the smoothing transform is given by any $\mu\in P
(\R)$ such that, if $R$ has distribution $\mu$, the equation
%
\begin{equation}
\label{equation} R=_d \sum_{i=1}^N
A_i R_i + B,
\end{equation}
holds true. We are going to distinguish between the case of $B=0$ a.s.
(the homogeneous smoothing transform) and the
other one called the nonhomogeneous smoothing transform.

The homogeneous equation \eqref{equation} is used for example, to study
interacting particle systems \cite{DL} or the branching random walk
\cite{HS,AR}.
In recent years, from practical reasons, the inhomogeneous
equation has gained importance. It appears for example, in the stochastic
analysis of the Pagerank algorithm (which is the heart of the
Google engine) \cite{JC1,JC2,VL} as well as in the analysis of a
large class of divide and conquer algorithms including the
Quicksort algorithm \cite{NR,R}. Both the homogeneous and the
inhomogeneous equation were recently used to describe equilibrium
distribution of a class of kinetic models and used for example, to
study the distribution of particle velocity in Maxwell gas (see, e.g.,
\cite{bassetti}).

Properties of the fixed points of equation \eqref{equation} are
governed by the function
\[
m(s) = \E \Biggl[ \sum_{i=1}^N
|A_i|^s \Biggr] = N \E\bigl[ |A_1|^s
\bigr].
\]

Suppose that $s_1 = \sup\{s\dvt  m(s)<\infty \}$ is strictly positive.
%
Clearly $m$ is convex and differentiable on $(0,s_1)$. We
assume that there are $0<\g<\a<s_1$ such that
\[
m(\g)=m (\a)=1.
\]
Then
\[
0<m'(\a)=\E \Biggl[\sum_{i=1}^N
|A_i|^{\a}\log|A_i| \Biggr]
\]
and the latter
quantity is finite. The main result of this paper is
the following theorem.

\begin{theorem}
\label{mthm} Suppose that
\begin{itemize}[$\bullet$]
\item[$\bullet$]$\log|A_1|$ is nonlattice;
\item[$\bullet$]$\P[A_1 > 0] > 0$ and $\P[A_1<0]>0$;
\item[$\bullet$]$s_1>0$;
\item[$\bullet$] there are $0<\g<\a<s_1$ such that $m(\g)=m (\a)=1$;
\item[$\bullet$] there is $\eps>0$ such that $\E|B|^{\g+\eps}<\infty $.
\end{itemize}
Suppose that $R$ is a nontrivial solution to \eqref{equation}
such that $\E|R|^{\g+\eps}<\infty $. Then
\[
\liminf_{t\to\infty } t^\a\P[R>t] >0 \quad \mbox{and}\quad
\liminf_{t\to\infty } t^\a\P[R<-t] >0.
\]
\end{theorem}

\begin{rem}
Under the assumptions of Theorem~\ref{mthm} the random variable $R$ is
real valued and it attains both positive and negative values. If $\P
[A_1>0] = \P[B>0]=1$ then $R$ is a positive random variable and
exactly the same proof shows that
\[
\liminf_{t\to\infty } t^\a\P[R>t] >0.
\]
\end{rem}

Existence of such a solution implies $\g<2$
for the nonhomogeneous case and $1\leq\g<2$ for the homogeneous
one (see \cite{ADM}). Then the solution is basically unique (given the
mean of it
exists) and, if $\E|B| ^{\a}<\infty $ then for every $s<\a$
%
\begin{equation}
\label{moment} \E|R|^s <\infty .
\end{equation}
In view of the result of Jelenkovic
and Olvera-Cravioto (Theorem~4.6 in \cite{JC3}), Theorem~\ref{mthm}
implies.

\begin{cor}
Suppose that the assumptions of Theorem~\ref{mthm} are satisfied
and additionally let $\E|B|^{\a}<\infty $.
Then
%
\begin{equation}
\label{asym} \lim_{t\to\infty } t^\a\P[R>t]=\lim
_{t\to\infty } t^\a\P [R<-t]=K >0.
\end{equation}
\end{cor}

The existence of the limit in \eqref{asym} for such $R$, in a more
general case
of random $N$, was proved by Jelenkovic and Olvera-Cravioto
\cite{JC3}, Theorem~4.6, but from the expression for $K$, given by
their renewal theorem, it is not possible to conclude its strict
positivity except of the very particular case when $A_1,\ldots, A_N,B$ are
positive and $\a
\geq1$. There are other solutions to \eqref{equation} than those
mentioned in the above corollary. For the full description of them
see \cite{ABM,AM,AM2}. Clearly, Theorem~\ref{mthm} matters only
for solutions satisfying \eqref{moment}.

Some partial results concerning positivity of $K$ are contained in
\cite{BDMM} and \cite{ADM}. The paper \cite{BDMM} deals with
matrices but Theorem~2.12 and Proposition~2.13 there can be
specified to our case. Under additional assumption that $\E
|B|^{s_0}<\infty $ for some $\a<s_0<s_1$ they say that either $K>0$ or
$\E|R|^s<\infty $ for all $s<s_0$. If $R$ is not constant, the latter
is not possible when there is $\b\leq s_0$ such that $\E
|A_1|^{\b}=1$. Indeed, then $R$ becomes the solution of
\[
R=AR+Q
\]
with $Q=\sum_{i=2}^NA_iR_i+B$ and the conclusion of
Goldie's theorem \cite{Go} would be violated. It is interesting that
for the
asymptotics in \eqref{asym} in the case of $N$ being constant the
implicit renewal theorem of Jelenkovic and Olvera-Cravioto is not
needed. The usual one on $R$ is sufficient \cite{BDMM}, Theorem~2.8. For positivity of $K$ in the general case of random $N$ see
\cite{ADM}, Theorem~9.

Clearly, Theorem~\ref{mthm} improves considerably the results of
\cite{BDMM} specialised to the one dimensional case. Also, the technique
is purely probabilistic while in \cite{BDMM} holomorphicity of $\E
|R|^z $ and the Landau theorem is used.


Let $\mu_A$ be the law of $A_i$. In Section~\ref{sec2}, we show some
necessary properties of the random walks with the transition
probability $\mu_A$. A version of the Bahadur, Rao theorem (\cite
{DZ}, Theorem
3.7.4) is needed and its proof is included in the \hyperref[appendix]{Appendix}.
Section~\ref{sec3} is devoted to the proof of Theorem~\ref{mthm}.

\section{Random walk generated by the measure \texorpdfstring{$\mu_A$}{$mu_A$}}\label{sec2}
In this section, we will study
properties of the random walk $\{|\wt A_1\cdots \wt A_n|\}_{n\in\N}$,
where $\wt A_i$ are independent and distributed according to the
measure $\mu_A$ (it is convenient for our purpose to use the
multiplicative notation). Since $\E\log|\wt A_1|<0$, by the
strong law of large numbers, this random walk converges to 0 a.s.
Nevertheless, our aim is to describe a sufficiently large set on
which trajectories of the process exceed an arbitrary large, but
fixed number $t$. Given $n$, one can prove that the probability of the
event $ \{|\wt A_1\cdots \wt A_n|>t \}$ is largest when $n$ is
comparable with $n_0$ defined by
%
\begin{equation}
\label{n0} n_0 = \biggl\lfloor\frac{\log t}{N\rho} \biggr\rfloor,
\end{equation}
where $\rho= \E [ |\wt A_1|^\a\log|\wt A_1|  ]$.
Notice that $n_0$ depends on $t$. However, since we are interested only
in estimates from below we need less and for our purpose it is
sufficient to consider
sets
%
\begin{equation}
\label{vn} {V}_{n} = \bigl\{ |\wt A_1\cdots \wt
A_n|\ge t \mbox{ and } |\wt A_1\cdots \wt A_s|\le
\mathrm{e}^{-(n-s)\d}tC_0 \mbox{ for every }s\le n-1 \bigr\},
\end{equation}
where $C_0$ is a large constant and $\d$ is a small constant (both
will be defined later).

Our main result of this section is the following
theorem.

\begin{theorem}
\label{thm_Vn} Assume $\E[|\wt A_1|^{\a+\d}]<\infty $, $\E[|\wt
A_1|^{\a}]=\frac{1}N$ and $0<\rho<\infty $.
There are constants $C_0, C_1, C_2$ such that for
sufficiently large $t$ and
for $n_0-\sqrt{n_0}\leq n \leq n_0$
\[
\frac{C_1}{\sqrt n t^\a N^n} < \P [ V_n ] \le \frac{C_2}{\sqrt n t^\a N^n}.
\]
\end{theorem}

In order to prove the theorem above we will need precise estimates
of $\P[|\wt A_1\cdots \wt A_s|>t]$. We will use the following extension
of the Bahadur, Rao theorem (\cite{DZ}, Theorem 3.7.4, see also Example
3.7.10).

\begin{prop}
\label{prop_br}
Assume $\E[|\wt A_1|^{\a+\d}]<\infty $, $\E[|\wt A_1|^{\a}]=\frac{1}N$ and $0<\rho<\infty $.
There is $C$ such that for every $d\geq0$ and every
$n\in\N$
%
\begin{equation}
\label{eq:p1} \P \bigl\{ |\wt A_1\cdots \wt A_n|>\mathrm{e}^d
\mathrm{e}^{ \rho nN} \bigr\}\leq \frac{C}{\sqrt{2\pi}\a\lambda\sqrt n \mathrm{e}^{\rho\a nN}N^{n}
\mathrm{e}^{\a d}} ,
\end{equation}
where $\lam=\sqrt{\Lambda''(\a)}$ for
$
\Lambda(s) = \log\E [|\wt A_1|^s ]$.

Moreover, let $\theta\geq0$ and 
\begin{equation}
\label{eq: p222} 0\leq\frac{d}{\sqrt n}\leq\theta
\end{equation}
for sufficiently large $n$.
Then there is $C=C(\theta)$ such that for large $n$:
%
\begin{equation}
\label{eq:p2} \sqrt{2\pi}\a\lambda\sqrt n \mathrm{e}^{\rho\a n N}N^n
\mathrm{e}^{\a d} \mathrm{e}^{\afrac{d^2}{2\lambda^2 n}}\cdot\P \bigl\{ |\wt A_1\cdots \wt
A_n|>\mathrm{e}^d \mathrm{e}^{\rho nN} \bigr\} = 1 + C(\theta)\mathrm{o}(1),
\end{equation}
where as usual $\lim_{n\to\infty } \mathrm{o}(1)=0$ uniformly for $d$
satisfying \eqref{eq: p222}.
\end{prop}

The proof is a slight modification of the proof of Theorem~3.7.4
in \cite{DZ}. For reader's convenience we give all the details but
we postpone the proof to the \hyperref[appendix]{Appendix}.

We will also use the following. Since $\E [|\wt A_1|^\b
]<\frac{1}N$ for $\b<\a$ and sufficiently close to $\a$, one can find
$\b<\a$ and $\g>0$ such that
%
\begin{equation}
\label{eq:2.4} \E \bigl[|\wt A_1|^\b \bigr] =
\frac{1}{N^{1+\g}}.
\end{equation}
\begin{pf*}{Proof of Theorem~\ref{thm_Vn}}
Denote
\begin{eqnarray*}
U_n &=& \bigl\{ |\wt A_1\cdots \wt A_n|>t \bigr\},
\\
W_{s,n} &=& \bigl\{ |\wt A_1\cdots \wt A_s| >
\mathrm{e}^{-\d(n-s)} C_0 t \bigr\}.
\end{eqnarray*}

We have
\begin{eqnarray*}
\P [ V_n ] &=& \P \biggl[ U_n \cap\bigcap
_{s<n} W_{s,n}^c \biggr]
\\
&=& \P [ U_n ] - \P \biggl[ U_n \cap \biggl(\bigcap
_{s<n} W_{s,n}^c
\biggr)^c \biggr]
\\
&=& \P [ U_n ] - \P \biggl[\bigcup_{s<n}
( U_n\cap W_{s,n} ) \biggr].
\end{eqnarray*}
By Proposition~\ref{prop_br} ($s=n, d = N\rho(n_0-n)$,
$\theta=N\rho+1$)
\begin{eqnarray*}
\P [ U_n ] &=& \P \bigl[ |\wt A_1\cdots \wt A_n|>t \bigr]
= \P \bigl[ |\wt A_1\cdots \wt A_n| > \mathrm{e}^{N\rho n}\mathrm{e}^{N\rho(n_0-n)}
\bigr]
\\
&\geq& \frac{C_1 \mathrm{e}^{-N\rho\a(n_0-n)}}{\sqrt n \mathrm{e}^{N\rho\a n}N^n} =
 \frac{C_1}{ \sqrt n \mathrm{e}^{N\rho\a n_0}N^n} = \frac{C_1}{\sqrt n t^\a
N^n}
\end{eqnarray*}
for sufficiently large $t$ and $C_1 = \frac{1+C(N\rho
+1)\mathrm{o}(1)} {\sqrt{2\pi}\a\lam}\exp ({-\frac{(N\rho
+1)^2}{2\lambda
^2}} )$. Exactly in the same way \eqref{eq:p2} gives estimates from
above with $C_2 = \frac{1+C(N\rho+1)\mathrm{o}(1)} {\sqrt{2\pi}\a\lam}$.
Therefore to prove the theorem, it is sufficient to justify that
%
\begin{equation}
\label{square} \P \biggl[\bigcup_{s<n} (
U_n\cap W_{s,n} ) \biggr] \le \frac{\eps}{\sqrt n t^\a N^n}.
\end{equation}

We fix $t$, $n_0$ and $n$ such that
$n_0-\sqrt{n_0}\le n \le n_0$.
First we estimate $\P[ U_n\cap W_{s,n} ]$ for $s< n-D\log n$, where
the constant $D$ will be defined later. By the Chebyshev inequality and
\eqref{eq:2.4}, we have
\begin{eqnarray*}
\P[ U_n\cap W_{s,n} ] &=& \sum
_{m=0}^\infty \P \bigl[\mathrm{e}^m
\mathrm{e}^{-\d
(n-s)}C_0 t < |\wt A_1\cdots \wt A_s|
\le \mathrm{e}^{m+1} \mathrm{e}^{-\d(n-s)}C_0 t \mbox{ and } |\wt
A_1\cdots \wt A_n| > t \bigr]
\\
&\le& \sum_{m=0}^\infty \P \bigl[ |\wt
A_1\cdots \wt A_s| > \mathrm{e}^{m} \mathrm{e}^{-\d
(n-s)}C_0
t \bigr] \P \bigl[ |\wt A_{s+1}\cdots \wt A_n| >
\mathrm{e}^{-(m+1)}\mathrm{e}^{\d(n-s)}C_0^{-1} \bigr]
\\
&\le& \sum_{m=0}^\infty \frac{ \mathrm{e}^{\d\a(n-s)} }{\mathrm{e}^{m\a}C_0^\a
t^\a}
\bigl(\E|\wt A_1|^\a \bigr)^s \cdot
\frac{\mathrm{e}^{\b(m+1)}C_0^\b
}{\mathrm{e}^{\d\b(n-s)}} \bigl(\E|\wt A_1|^\b
\bigr)^{n-s}
\\
&\le& \frac{\mathrm{e}^{\d\a(n-s)}}{C_0^{\a-\b}t^\a} \cdot\frac{1}{N^s} \cdot\frac{1}{ \mathrm{e}^{\d\b(n-s)} N^{n-s} N^{\g(n-s)} } \cdot
\sum_{m=0}^\infty \frac{\mathrm{e}^\b}{\mathrm{e}^{m(\a-\b)}}
\\
&\le& \frac{C}{C_0^{\a-\b}t^\a N^n \mathrm{e}^{(\g\log N+\d(\b-\a
))(n-s)}} = \frac{C}{C_0^{\a-\b}t^\a N^n \mathrm{e}^{\g_1(n-s)}},
\end{eqnarray*}
where
$
\g_1:=\g\log N + (\b-\a)\d
$ and choosing appropriately small $\d$ we can assume that $\g_1>0$.
Hence, for $s<n-D\log n$
%
\begin{equation}
\label{star} \P[ U_n\cap W_{s,n} ] \le
\frac{C}{C_0^{\a-\b}t^\a N^n \mathrm{e}^{\g_1(n-s)}}.
\end{equation}

For $s> n-D\log n$, we estimate
\begin{eqnarray*}
\P[ U_n\cap W_{s,n} ] &=& \sum_{m=0}^\infty
\P \bigl[\mathrm{e}^m \mathrm{e}^{-\d
(n-s)}C_0 t < |\wt
A_1\cdots \wt A_s| \le \mathrm{e}^{m+1} \mathrm{e}^{-\d(n-s)}C_0
t \mbox{ and } |a_1\cdots a_n| > t \bigr]
\\
&\le&\sum_{m=0}^\infty \P \bigl[ |\wt
A_1\cdots \wt A_s| > \mathrm{e}^m \mathrm{e}^{-\d
(n-s)}C_0
t \bigr] \P \bigl[ |\wt A_{s+1}\cdots \wt A_n| >
\mathrm{e}^{-(m+1)}\mathrm{e}^{\d(n-s)}C_0^{-1} \bigr].
\end{eqnarray*}
We denote the first factor of the sum by $I_m$. To estimate it, we will
use Proposition~\ref{prop_br}. Namely let
\begin{eqnarray*}
k&=& n-s, \qquad k_0 = n_0 - s,
\\
d_1 &=& -\d k + m + \log C_0 + N\rho
k_0,
\\
d_2 &=& d_1 + 1,
\end{eqnarray*}
then (recall $\log t = (s+k_0)N\rho$)
\[
\mathrm{e}^m \mathrm{e}^{-\d(n-s)}C_0 t = \mathrm{e}^{d_1}
\mathrm{e}^{N\rho s}.
\]
So, by Proposition~\ref{prop_br}:
\[
\P \bigl[ |\wt A_1\cdots \wt A_s| > \mathrm{e}^{d_1+N\rho s} \bigr]
\le\frac
{C}{\sqrt s} N^{-s}\mathrm{e}^{-N\rho\a s-\a d_1} \leq \frac{C\mathrm{e}^{\d\a k}}{ C_0^\a
\mathrm{e}^{\a m} t^\a N^s \sqrt s}.
\]
The second factor we estimate exactly in the same way as previously and
we obtain
%
\begin{eqnarray}
\label{2star} %
\P[ U_n\cap W_{s,n} ] &=& \sum
_{m=0}^\infty \frac{C\mathrm{e}^{\d\a(n-s)}}{
C_0^\a \mathrm{e}^{\a m} t^\a N^s \sqrt s} \cdot
\frac{ \mathrm{e}^{\b(m+1)} C_0^\b}{ \mathrm{e}^{\d\b(n-s)} N^{(1+\g)(n-s)}
}\nonumber
\\[-8pt]\\[-8pt]
&\le&\frac{C}{ C_0^{\a-\b} t^\a N^n \sqrt n \mathrm{e}^{\g_1(n-s)} }. \nonumber%
\end{eqnarray}
%

Next, in view of \eqref{star} and \eqref{2star}
\begin{eqnarray*}
\P \biggl[\bigcup_{s<n} ( U_n\cap
W_{s,n} ) \biggr] &\le& \sum_{s< n-D\log n} \P [
U_n\cap W_{s,n} ] + \sum_{ n-D\log n\le s < n}
\P [ U_n\cap W_{s,n} ]
\\
&\le& \sum_{s< n-D\log n} \frac{C}{C_0^{\a-\b}t^\a N^n \mathrm{e}^{\g_1(n-s)}} + \sum
_{ n-D\log n\le s < n} \frac{C}{C_0^{\a-\b}t^\a N^n \sqrt n
\mathrm{e}^{\g_1(n-s)}}
\\
&\le& \frac{C}{C_0^{\a-\b}t^\a N^n} \biggl( \frac{n}{n^{\g_1 D }} +\frac{1}{\sqrt n} \sum
_{ s < D\log n} \frac{1}{\mathrm{e}^{\g_1 s}} \biggr)
\\
&\le& \frac{C}{C_0^{\a-\b}t^\a N^n} \biggl( \frac{1}{n^{\g_1 D -1}} +\frac{1}{\sqrt n} \biggr)
\\
&\le& \frac{\eps}{\sqrt n t^\a N^n }
\end{eqnarray*}
assuming that $\frac{C}{C_0^{\a-\b}}<\eps$ and $\g_1 D \ge\frac{3}2$. Hence, \eqref{square} and the proof is finished.
\end{pf*}

\section{Proof of Theorem \texorpdfstring{\protect\ref{mthm}}{1.1}}\label{sec3}

We start with the following lemma.

\begin{lem} If $\E [|A_1|^\b\log|A_1| ] >0$ for some $\b
>0$, $\P[A_1>0]>0$ and $\P[A_1<0]>0$, then any nontrivial solution of
\eqref{equation} is unbounded at $+\infty $ and $-\infty $.
\end{lem}

\begin{pf}
Suppose that $R$ is a bounded solution of \eqref{equation} and $R\neq C$ a.s. for any $C$. Assume first that $R$ is bounded a.s. from below
and from above. Let $[r,s]$ be the smallest interval containing the
support of $R$ for some finite numbers $r$ and $s$. Of
course $r\neq s$.
Denote $\wt B = \sum_{i=2}^N A_i R_i+B$, then
%
\begin{equation}
\label{eqaffine} R =_d A_1R_1+\wt B.
\end{equation}
Since $\E[|A_1|^\b\log|A_1|]>0$, the probability of the set $ U =
\{ (A_1,\wt B)\dvt  |A_1| > 1  \} $ is strictly positive.
Then by \eqref{eqaffine}
we must have
\[
A_1 r +\wt B \ge r \quad \mbox{and}\quad  A_1 s +\wt B \le s\qquad
\mbox{a.s.}
\]
But if we take a random pair $(A_1,\wt B)\in U$, then
\[
\bigl|(A_1 r + \wt B) - (A_1 s + \wt B)\bigr| =
|A_1||r-s| > |r-s|.
\]
Thus, we are led to a contradiction and at least one constant $r$ or
$s$ must be infinite. Without loss of generality, we can assume that
$s=+\infty $. In view of our assumptions, we can choose a large
constant $M$ and a small constant $\eps$ such that the probability of
the set $ V =
\{ (A_1,\wt B)\dvt  A_1<-\eps, \wt B < M  \} $ is strictly
positive. Now, take any $x> (r-M)/(-\eps)$ belonging to the support of
$R$. Then for any $(A_1,\wt B)\in V$ we have
\[
A_1 x+\wt B < -\eps x+M < r.
\]
Thus, by \eqref{eqaffine}, $r$ cannot be a lower bound of the support
of $R$ and must be equal $-\infty $.
%
\end{pf}

Let $\T$ be an $N$-ary rooted tree, that is, the tree with a
distinguished vertex $o$ called root, such that every vertex has $N$
daughters and one mother
(except the root). The tree $\T$ can be identified with the set of
finite words over the alphabet $\{1,2,\ldots ,N\}$:
\[
\T= \bigcup_{k=0}^\infty \{1,2,\ldots ,N
\}^k,
\]
where the empty word $\emptyset$ is the root and given $i_1i_2\cdots i_n\in
\T$ its daughters are the words of the form $i_1i_2\cdots i_nj$ for $j=1,\ldots ,N$.
We denote a typical vertex of the tree by $\g= i_1i_2\cdots i_n$ and we
identify it with the shortest path connecting $\g$ with $o$. We write
$|\g|=n$
for the length of $\g$ and $\g_{|_k} = i_1\cdots i_k$ for the curtailment
of $\g$ after $k$ steps. Conventionally, $|\emptyset|=0$ and $\g
_{|_0} = \emptyset$. If $\g_1 = i_1^1i_2^1\cdots i_{n_1}^1\in\T$ and $\g
_2 = i_1^2i_2^2\cdots i_{n_2}^2\in\T$ then we write $\g_1\g_2 =
i_1^1i_2^1\cdots i_{n_1}^1i_1^2i_2^2\cdots i_{n_2}^2 $
for the element of $\T$ obtained by juxtaposition. In particular, $\g
\emptyset= \emptyset\g= \g$. We partially order $\T$ by writing
$\g_1\le\g_2$ if there exists $\g_0\in\T$ such that $\g_2 = \g
_1\g_0$. For two vertices $\g_1$ and $\g_2$, we denote by $\g_0 =
\g_1\wedge\g_2$ the longest common subsequence of $\g_1$ and $\g
_2$ that is, the maximal $\g_0$ such that both $\g_0\le\g_1$ and
$\g_0\le\g_2$.

To every vertex $\g\in\T$ we associate random variables $(A_{\g
1},\ldots ,A_{\g N},B_\g,R_{\g1},\ldots ,R_{\g N})$ which are independent
copies of
$(A_{1},\ldots ,A_{N},B,R_{1},\ldots ,R_{N})$ defined in \eqref{equation}. It is
more convenient to think that $A_{\g i}$ and $R_{\g i}$ are indeed
attached not to the vertex $\g$ but to the edge connecting $\g$ with
$\g i$.
We write $\Pi_\g= A_{\g_{|_1}}A_{\g_{|_2}}\cdots A_{\g}$, then $\Pi_\g
$ is just the product of random variables $A_{\g_{|_k}}$ which are
associated with consecutive edges connecting the root $o$ with $\g$.

We fix $\g= i_1\cdots i_n$ and we apply $n$ times the stochastic equation
\eqref{equation} in such a way that in $k$th step we apply recursively
this equation to $R_{\g_{|_k}}$:
%
\begin{eqnarray}
\label{eq:2star} %
R &=_d& \sum_{i=1}^N
A_i R_i + B_0\nonumber
\\
&=_d& A_{i_1} \Biggl( \sum_{j=1}^N
A_{i_1j}R_{i_1j} + B_{i_1} \Biggr) + \sum
_{i\neq i_1} A_i R_i + B_0\nonumber
\\
&=_d& A_{i_1}A_{i_1i_2}R_{i_1i_2} +
A_{i_1} \Biggl( \sum_{j\neq i_2}^N
A_{i_1j}R_{i_1j} + B_{i_1} \Biggr) + \sum
_{i\neq i_1} A_i R_i + B_0\nonumber
\\
&=_d& \Pi_{\g_{|_2}} R_{\g_{|_2}} + \sum
_{j\neq i_2}^N \Pi_{(\g
_{|_1}j)}R_{(\g_{|_1}j)} +
A_{i_1}B_{i_1} + \sum_{i\neq i_1}
A_i R_i + B_0
\\
&=_d& \Pi_{\g_{|_2}} \Biggl( \sum_{i=1}^N
A_{(\g_{|_2}i)}R_{(\g
_{|_2}i)} + B_{\g_{|_2}} \Biggr) + \sum
_{i\neq i_2}^N \Pi_{(\g
_{|_1}i)}R_{(\g_{|_1}i)} +
\sum_{i\neq i_1} A_i R_i +
A_{i_1}B_{i_1} + B_0\nonumber
\\
&=_d& \cdots\nonumber
\\
&=_d& \Pi_\g R_\g+ \sum
_{k<n} \sum_{i\neq i_k}
\Pi_{(\g
_{|_k}i)}R_{(\g_{|_k}i)} + \sum_{k<n}
\Pi_{\g_{|_k}}B_{\g_{|_k}}.\nonumber %
\end{eqnarray}
We define
\[
V_{\g} = \bigl\{ |\Pi_\g |\ge t \mbox{ and } |\Pi
_{\g_{|_s}} |\le \mathrm{e}^{-(|\g|-s)\d}C_0 t \mbox{ for every }s<|
\g| \bigr\}.
\]
Notice that if we denote $\wt A_k = A_{\g_{|_k}}$, then the set $V_\g
$ coincides with the set $V_{|\g|}$ defined in \eqref{vn}. Thus, by
Theorem~\ref{thm_Vn} we can choose large $C_0$ such that if $n=|\g|$
and $n_0 - \sqrt{n_0} < n < n_0$, then
\[
\P[V_\g] \ge\frac{C}{\sqrt n t^\a N^n}.
\]
For a sufficiently large constant $d$ (defined later) and
$D=\frac{Nd^2 + d}{1-\mathrm{e}^{-\sfrac{\d}{2}}}$, we define sets
\begin{eqnarray*}
W_\g&=& \bigl\{ |R_{(\g_{|_s}i)} | < d \mathrm{e}^{{(|\g|-s)\d
}/{4}},
|A_{(\g_{|_s}i)} | < d \mathrm{e}^{{(|\g|-s)\d}/{4}}, |B_{\g
_{|_s}}| < d
\mathrm{e}^{{(|\g|-s)\d}/{2}}, \\
&&\hphantom{\bigl\{}s=0,\dots,|\g|-1; i\neq i_{s+1} \bigr\};
\\
W^+_\g&=& W_\g\cap \{ R_\g> 2D \};
\\
W^-_\g&=& W_\g\cap \{ R_\g< -2D \};
\\
V^+_\g&=& V_\g\cap \{ \Pi_\g>0 \};
\\
V^-_\g&=& V_\g\cap \{ \Pi_\g<0 \}.
\end{eqnarray*}
Finally we define
\[
\wt V_\g= \bigl(V_\g^+ \cap W_\g^+ \bigr)
\cup \bigl(V_\g^- \cap W_\g^- \bigr).
\]

\begin{lem}
Assume $\g\in\T$. Then on the set $\wt V_\g$ we have
\[
R > At.
\]
\end{lem}

\begin{pf}
Let $n=|\g|$, then by \eqref{eq:2star} on $\wt V_{\g}$ we have
\begin{eqnarray*}
R &\ge& \Pi_\g R_\g- \biggl| \sum
_{k<n}\sum_{i\neq i_k} \Pi
_{(\g_{|_k}i)} R_{(\g_{|_k}i)} + \sum_{k<n}
\Pi_{\g_{|_k}} B_{\g_{|_k}} \biggr|
\\
&\ge& 2Dt - \sum_{k<n} \bigl(Nd^2+d
\bigr) \mathrm{e}^{-{(n-k)\d}/{2}}C_0 t
\\
&\ge& Dt.
\end{eqnarray*}
\upqed
\end{pf}
We are going to prove that for some $\eta>0$
%
\begin{equation}
\label{eq:3star} \P \biggl[\bigcup_{\{\g\in\T: n_0-\sqrt{n_0} < |\g| <n_0\}} \wt
V_\g \biggr]\ge\eta t^{-\a},
\end{equation}
which immediately
implies that
\[
\liminf_{t \to\infty }\P\{ R>t\}t^{\a}>0.
\]

\begin{lem}
Let $X_i$ be a sequence of i.i.d. random variables such that $\E
|X_1|^{\eps}<\infty $ for some $\eps>0$. Let $\d_0>0$. Then there
exist constants $d_0$ and $p_0>0$
such that for every $n$
\[
\P \bigl[ |X_i| < d_0 \mathrm{e}^{(n-i)\d_0}, i=1,2,\ldots ,n-1
\bigr] \ge p_0.
\]
\end{lem}

\begin{pf}
By the Chebyshev inequality, we have
\[
\P \bigl[ |X_i| \ge d_0 \mathrm{e}^{(n-i)\d_0} \bigr] \le
\frac{\E|X_i|^{\eps
}}{d_0^{\eps}} \mathrm{e}^{-(n-i)\d_0\eps}.
\]
Take $d_0$ such that $d_0^{\eps} > 3 \E|X_i|^{\eps}$. Then,
since $1-\frac{x}{3} > \mathrm{e}^{-x}$ for $x\in[0,1]$ we have
\[
\P \bigl[ |X_i| < d_0 \mathrm{e}^{(n-i)\d_0} \bigr] \ge1 -
\tfrac{1}{3} \mathrm{e}^{-(n-i)\d_0 \eps} \ge\exp \bigl({- \bigl(\mathrm{e}^{-\d_0 \eps}
\bigr)^{n-i}} \bigr).
\]
Therefore,
\begin{eqnarray*}
&&\P \bigl[ |X_i| < d_0 \mathrm{e}^{(n-i)\d_0}, i=1,2,\ldots ,n-1
\bigr] \\
&&\quad = \prod_{i=1}^{n-1} \P \bigl[
|X_i| < d_0 \mathrm{e}^{(n-i)\d_0} \bigr]\ge\prod
_{i=1}^{n-1} \mathrm{e}^{- (\mathrm{e}^{-\d_0 \eps})^{n-i}}
\\
&&\quad = \exp \Biggl({-\sum_{i=1}^{n-1}
\bigl(\mathrm{e}^{-\d_0\eps}\bigr)^i} \Biggr)\ge\exp \bigl({-
\bigl(1-\mathrm{e}^{-\d_0 \eps}\bigr)^{-1}} \bigr) =:p_0.
\end{eqnarray*}
\upqed
\end{pf}
Since $B$ and $R$ have absolute moments of order bigger then $\g$ we
obtain the following corollary.

\begin{cor}
There are constants $d$ and $p>0$ such that for every $\g\in\T$
\[
\P\bigl[W_\g^+\bigr] \ge p \quad \mbox{and}\quad  \P\bigl[W_\g^-
\bigr]\ge p.
\]
\end{cor}

In view of the last result to obtain \eqref{eq:3star}, it is
sufficient to prove
\[
\P \biggl[\bigcup_{\{\g\in\T: n_0-\sqrt{n_0} < |\g| <n_0\}} V_\g \biggr]
\ge\eta_1 t^{-\a},
\]
for some $\eta_1>0$.

In fact, we will estimate from below much smaller sum over a sparse
subset of $\T$. The details are as follows.

We fix a large integer $C_1$ (determined later) and an arbitrary
element $\ov\g$ of $\T$ such that $|\ov\g| = C_1$ (e.g., $\ov\g$
can be chosen as the word consisting of $n$ one's). We define a sparse
subset of vertices of $\T$:
\[
\ov\T= \bigl\{ \g\in\T\dvt  \bigl(|\g|\quad \mathrm{mod}\ C_1\bigr) = 0, \g= \g_{|_{|\g
|-C_1}}
\ov\g, n_0 - \sqrt{n_0} < |\g| < n_0
\bigr\},
\]
that is, $\ov\T$ is the set of vertices of $\T$ located on the level
$kC_1$ (for some integer $k$) such that
$n_0-\sqrt{n_0} < kC_1 < n_0$ and such that the last $n$ letters of
$\g$ form the word $\ov\g$. Notice that for every $\g$ such that
$|\g|=kC_1$ the set
\[
\Biggl\{ \g\g_1, \g_1\in\bigcup
_{i=1}^{C_1}\{1,\ldots ,N\}^i \Biggr\}
\]
contains exactly one element of $\T$. Thus there are exactly
$N^{kC_1}$ elements of $\ov\T$ of length $(k+1)C_1$. Moreover, the
crucial property of the set $\ov\T$, that will be strongly used
below, is that the distance between two different elements of $\ov\T$
is at least $C_1$ (by ``distance'' we mean the usual distance on graphs,
that is, the minimal number of edges connecting two vertices).

\begin{theorem} There is $\eta>0$ such that
\[
\P \biggl(\bigcup_{\gamma\in \ov\T} V_\g \biggr)
\geq\frac{C\eta
}{N^{C_1}C_1t^{\a}}.
\]
\end{theorem}

\begin{pf}
By the inclusion--exclusion principle, we have
%
\begin{equation}
\label{eq:1} \P \biggl(\bigcup_{\gamma\in\ov\T}
V_\g \biggr) \ge \sum_{\g\in\ov\T}\P (
V_\g ) - \sum_{\g\in\ov\T
}\sum
_{U_\g} \P ( V_\g\cap V_{\g'} ),
\end{equation}
where $U_\g= \{ \g'\in\ov\T\setminus\{\g\}\dvt  |\g'|\le|\g| \}$.

Therefore, we have to estimate
\[
\sum_{\g\in\ov\T}\P ( V_\g ) \quad \mbox{and}\quad
\sum_{\g\in
\ov\T}\sum_{ U_\g} \P
( V_\g\cap V_{\g'} ).
\]
Let $K$ be the set of levels on which there are some elements of $\ov
\T$, that is,
\[
K= \{ kC_1\dvt  n_0 - \sqrt{n_0} <
kC_1 < n_0 \}.
\]
Let $L=|K|$ be the number of elements of the set $K$ and let $n_j$ be
the $j$th element of $K$.

Observe that for given $n\in K$ there are
exactly $N^{n-C_1}$ elements of $\ov\T$ located on the level $n$ and
for every such $\g$, by Theorem~\ref{thm_Vn}, we have
$\P( V_\g)\geq
\frac{C}{\sqrt{n}N^nt^{\a}}$. Hence,
%
\begin{equation}
\label{eq:2} \sum_{\g\in\ov\T}\P(V_\g)\geq
\sum_{j=1}^L \frac{C}{\sqrt{n_j}N^{n_j}t^{\a}}N^{n_j-C_1}
\geq \frac{C }{N^{C_1}C_1t^{\a}}.
\end{equation}
Now, let us estimate the sum of intersections. We fix $\g\in\ov\T$
and $\g'\in U_\g$. Let $\g_0 = \g\wedge\g'$ and let $s$ be the
length of $\g_0$. We have
%
\begin{eqnarray}
\label{int} %
\P [ V_\g\cap V_{\g'} ] &\le&\P
\bigl[ V_\g\cap \bigl\{ |\Pi_{\g_0}| < \mathrm{e}^{-\d(|\g|-s)}C_0
t, |\Pi_{\g'}| > t \bigr\} \bigr]\nonumber
\\
&\le&\P [ V_{\g} ]\P \bigl[ | A_{\g'_{|_{s+1}}}A_{\g
'_{|_{s+2}}}\cdots A_{\g'}
| > \mathrm{e}^{\d(|\g|-s)} C_0^{-1} \bigr]
\\
&\le&\P [ V_{\g} ]\cdot\frac{C_0^\a}{\mathrm{e}^{\a\d(|\g|-s)}
N^{|\g'|-s}},\nonumber %
\end{eqnarray}
where for the last inequality we have used the Chebyshev
inequality. We fix $\g\in\ov\T$ and we consider $\g'\in U_\g$.
Notice that if $\g$ and $\g'$ connect on the level $s$, that is, $\g
_{|_s} = \g\wedge\g'$, then $s$ must be smaller than
$|\g|-C_1$. Given $s$ let us estimate the number of elements
$\g'\in U_\g$ such that $\g_{|_s} = \g\wedge\g'$. All these
elements must be located on levels $|\g|, |\g|-C_1,\ldots , |\g|-k
C_1$, where $k$ is the largest number such that $|\g|-kC_1 \ge
\max\{ s,n_0-\sqrt{n_0}\}$, that is,
\[
k \le\frac{1}{C_1} \min \bigl\{ |\g|-s , |\g|-n_0 +
\sqrt{n_0} \bigr\} \le\frac{1}{C_1} \bigl(|\g|-s\bigr).
\]
Moreover on the level $|\g|-jC_1$ ($j<k$), there are exactly $N^{|\g
|-jC_1-s-C_1}$ elements of $U_\g$. Thus for $C_1$ sufficiently large,
by \eqref{int}, we have
\begin{eqnarray*}
&&\sum_{\g\in\ov\T}\sum_{\g' \in U_\g}
\P [ V_\g\cap V_{\g
'} ] \\
&&\quad \le \sum
_{\g\in\ov\T} \sum_{s\le|\g|-C_1} \sum
_{\{
\g'\in U_\g: \g_{|_s} = \g\wedge\g'\}} \P[V_\g] \cdot\frac{C_0^\a}{\mathrm{e}^{\a\d(|\g|-s)} N^{|\g'|-s}}
\\
&&\quad \le \sum_{\g\in\ov\T} \P[V_\g] \sum
_{s\le|\g|-C_1} \sum_{0\le j\le\sklfrac{1}{C_1}(|\g|-s)} \sum
_{\{\g'\in U_\g: \g_{|_s} = \g\wedge\g', |\g'| = |\g
|-jC_1\}} \frac{C_0^\a}{\mathrm{e}^{\a\d(|\g|-s)} N^{|\g'|-s}}
\\
&&\quad \le \sum_{\g\in\ov\T} \P[V_\g] \sum
_{s\le|\g|-C_1} \sum_{0\le j\le\sklfrac{1}{C_1}(|\g|-s)}
\frac{C_0^\a}{\mathrm{e}^{\a\d(|\g|-s)} N^{|\g|-jC_1-s}} \cdot N^{|\g| -
jC_1 -s - C_1}
\\
&&\quad \le \sum_{\g\in\ov\T} \P[V_\g] \sum
_{s\le|\g|-C_1} \frac{C_0^\a(|\g|-s)}{ C_1 N^{C_1} \mathrm{e}^{\a\d(|\g|-s)} }
\\
&&\quad \le \sum_{\g\in\ov\T} \P[V_\g]
\frac{C_0^\a}{ C_1 N^{C_1} \mathrm{e}^{{\a\d C_1}/2 }} \le\frac{1}2 \sum_{\g
\in\ov\T}
\P[V_\g].
\end{eqnarray*}
Finally, combining the above estimates with \eqref{eq:1} and \eqref
{eq:2}, we obtain
\[
\P \biggl[ \bigcup_{\g\in\ov\T} V_\g \biggr]
\ge\frac{1}2 \frac
{C}{N^{C_1} C_1 t^\a}.
\]
\upqed
\end{pf}

\begin{appendix}
\section*{Appendix: Proof of Proposition \texorpdfstring{\lowercase{\protect\ref{prop_br}}}{2.2}}
\label{appendix}
\setcounter{equation}{0}
\begin{pf}
We proceed as in \cite{DZ} and for reader's convenience we use the
same notation. Define
$X_i = \log|\wt A_i|$ and $\wh S_n = \frac{1}n\sum_{i=1}^n X_i$.
We introduce a new probability measure: $\wt\mu(dx) = N \mathrm{e}^{\a x
}\mu(dx)$, where $\mu$ is the law of $X_i$. Next, we normalize
$X_i$ and we define
new random variables: $Y_i =
\frac{X_i-N\rho}{\sqrt{\Lambda''(\a)}}$ and $W_n= \frac{1}{\sqrt
n}\sum_{i=1}^n Y_i$. Then $\E_{\wt\mu} Y_i=0$
and
\[
\wh S_n - N\rho=\frac{\lam}{\sqrt{n}}\frac{1}{\sqrt{n}}\sum
_{i=1}^n Y_i= \frac{\lam}{\sqrt{n}}W_n,
\]
where $\lam= \sqrt{\Lambda''(\a)}$ and $\Lambda(s) = \log (
\E [ |\wt A_1|^s  ]  )$.
Let $F_n$
be the distribution of $W_n$ with respect to the changed measure $\wt
\mu$.
Let $\psi_n = \a\lambda\sqrt n$.
Then,
\begin{eqnarray*}
\P \bigl\{ |\wt A_1\cdots \wt A_n|>\mathrm{e}^d
\mathrm{e}^{N\rho n} \bigr\} &=& \P \{ \wh S_n > N\rho+ d/n \}
\\
&=& \P \biggl\{ W_n> \frac{d}{\lam\sqrt {n}} \biggr\} =\E_{\wt
\mu}
\bigl[ N^{-n} |\wt A_1\cdots\wt A_n|^{-\a}
\mathbf{1}_{\{W_n >
\afrac{d}{\lambda
\sqrt n}\}} \bigr]
\\
&=& \mathrm{e}^{-\a n\rho N} N^{-n} \E_{\wt\mu} \bigl[\mathrm{ e}^{-\psi_n W_n}
\mathbf{1}_{\{W_n > \afrac{d}{\lambda\sqrt n}\}} \bigr]
\\
&=& \mathrm{e}^{-\a n\rho N} N^{-n} \int_{\afrac{d}{\lambda\sqrt
n}}^\infty
\mathrm{e}^{-\psi_n x}\,\mathrm{d}F_n(x).
\end{eqnarray*}
We will use here the Berry--Esseen expansion for nonlattice
distributions of $F_n$ (see \cite{F}, page 538):
%
\begin{equation}
\label{Berry-Essen} \lim_{n\to\infty } \biggl( \sqrt n \sup
_x \biggl| F_n(x) - \Phi(x) - \frac{m_3}{6\sqrt n}
\bigl(1-x^2\bigr)\phi(x) \biggr| \biggr) = 0,
\end{equation}
where $m_3 = \E_{\wt\mu}[Y_1^3]<\infty $, $\phi(x)=\frac{1}{\sqrt {2\pi}} \mathrm{e}^{-\sfrac{x^2}2}$ is the standard normal density, and $\Phi
(x)=\int_{-\infty }^x \phi(y)\,\mathrm{d}y$ is its distribution function.

First, we integrate by parts and then we use the above result
\begin{eqnarray*}
J&=&\a\lambda\sqrt n \mathrm{e}^{N\rho\a n}N^n \P \bigl\{ |\wt
A_1\cdots \wt A_n|>\mathrm{e}^d \mathrm{e}^{N\rho n} \bigr\}
\\
&=& \int_{\afrac{d}{\lambda\sqrt n}}^\infty
\psi_n \mathrm{e}^{-\psi_n
x}\,\mathrm{d}F_n(x)
\\
&=& \psi_n \mathrm{e}^{-\psi_n x} F_n(x) \biggl|_{\afrac{d}{\lambda\sqrt
n}}^\infty
+ \int_{\afrac{d}{\lambda\sqrt n}}^\infty \psi _n^2\mathrm{e}^{-\psi_n x}F_n(x)\,\mathrm{d}x
\\
&=& - \psi_n \mathrm{e}^{-\a d} F_n \biggl(
\frac{d}{\lambda\sqrt n} \biggr) + \int_{\a d}^\infty
\psi_n \mathrm{e}^{- x}F_n \biggl(\frac{x}{\psi_n}
\biggr)\,\mathrm{d}x
\\
&=& \int_{\a d}^\infty \psi_n
\mathrm{e}^{- x} \biggl[ F_n \biggl(\frac{x}{\psi
_n} \biggr) -
F_n \biggl(\frac{d}{\lambda\sqrt n} \biggr) \biggr]\, \mathrm{d}x
\\
&=& \mathrm{o}(1) \mathrm{e}^{-\a d} + \int_{\a d}^\infty
\psi_n \mathrm{e}^{- x} \biggl[ \Phi \biggl(\frac{x}{\psi_n}
\biggr) - \Phi \biggl(\frac{d}{\lambda\sqrt n} \biggr) \biggr] \,\mathrm{d}x
\\
&&{} + \frac{m_3 }{6\sqrt n} \int_{\a d}^\infty
\psi_n \mathrm{e}^{- x} \biggl[ \biggl(1- \biggl(\frac{x}{\psi_n}
\biggr)^2 \biggr)\phi \biggl(\frac{x}{\psi
_n} \biggr) - \biggl(1-
\biggl(\frac{d}{\lambda\sqrt n} \biggr)^2 \biggr) \phi \biggl(
\frac{d}{\lambda\sqrt n} \biggr) \biggr] \,\mathrm{d}x.
\end{eqnarray*}
We denote the second term by $I(n)$ and the third one by $\mathit{II}(n)$. Thus,
\[
J(n) = \mathrm{o}(1)\mathrm{e}^{-\a d} + I(n) + \mathit{II}(n).
\]
We estimate first $I$:
\begin{eqnarray*}
\sqrt{2\pi} I(n) &=& \int_{\a d}^\infty
\psi_n \mathrm{e}^{-x}\int_{\afrac{d}{\lambda\sqrt n}}^{\afrac{x}{\psi_n}}
\mathrm{e}^{-\sfrac{y^2}{2}}\,\mathrm{d}y\,\mathrm{d}x = \int_{\afrac{d}{\lambda\sqrt n}}^\infty
\psi_n \mathrm{e}^{-\sfrac
{y^2}{2}}\int_{ \psi_n y}^{\infty }
\mathrm{e}^{-x}\,\mathrm{d}x\,\mathrm{d}y
\\
&=& \int_{\afrac{d}{\lambda\sqrt n}}^\infty \psi_n
\mathrm{e}^{-\psi_n
y}\mathrm{e}^{-\sfrac{y^2}{2}}\,\mathrm{d}y = - \mathrm{e}^{-\psi_n y} \mathrm{e}^{-\sfrac{y^2}{2}}
\biggl|_{\afrac{d}{\lambda\sqrt
n}}^\infty - \int_{\afrac{d}{\lambda\sqrt n}}^\infty
y \mathrm{e}^{-\psi_n
y}\mathrm{e}^{-\sfrac{y^2}{2}}\,\mathrm{d}y
\\
&=& \mathrm{e}^{-\a d} \mathrm{e}^{-\afrac{d^2}{2\lambda^2 n}} - \int_{\afrac{d}{\lambda
\sqrt n}}^\infty
y \mathrm{e}^{-\psi_n y}\mathrm{e}^{-\sfrac{y^2}{2}}\,\mathrm{d}y.
\end{eqnarray*}
Let $\d>0$. We divide the last integral into two parts
\[
\int_{\afrac{d}{\lambda\sqrt n}}^\infty y \mathrm{e}^{-\psi_n y}\mathrm{e}^{-\sfrac
{y^2}2}\,\mathrm{d}y
= \int_{\afrac{d}{\lambda\sqrt n}}^{{\afrac{d}{\lambda\sqrt n}} + \sfrac\d
{\lambda}} y \mathrm{e}^{-\psi_n y}\mathrm{e}^{-\sfrac{y^2}{2}}\,\mathrm{d}y
+ \int_{{\afrac{d}{\lambda\sqrt n}} + \sfrac\d{\lambda}}^\infty y \mathrm{e}^{-\psi_n y}\mathrm{e}^{-\sfrac{y^2}{2}}\,\mathrm{d}y
\]
and denote the first one by $I_1(n)$ and the second one by $I_2(n)$. Then
\[
\mathrm{e}^{\a d} I_1(n) = \int_{\afrac{d}{\lambda\sqrt n}}^{\afrac{d}{\lambda
\sqrt n}+ \sfrac\d{\lambda}}
y \mathrm{e}^{\a d-\psi_n y}\mathrm{e}^{-\sfrac
{y^2}{2}}\,\mathrm{d}y \le\frac{\d}{\lambda}\frac{\theta+{\d}}{\lambda}
\cdot \mathrm{e}^{-\afrac{d^2}{2\lambda^2 n}}
\]
and large $n$ we have
\[
\mathrm{e}^{\a d}I_2(n) = \int_{\afrac{d}{\lambda\sqrt n}+ \sfrac\d{\lambda
}}^\infty
y \mathrm{e}^{\a d-\psi_n y}\mathrm{e}^{-\sfrac{y^2}{2}}\,\mathrm{d}y \le \mathrm{e}^{-\a\d\sqrt n} \mathrm{e}^{-\afrac{d^2}{2\lambda^2 n}}
\le\d \mathrm{e}^{-\afrac{d^2}{2\lambda^2 n}}.
\]
Thus, we have proved that for large $n$
\[
\sqrt{2\pi} \mathrm{e}^{\a d} I(n) = \mathrm{e}^{-\afrac{d^2}{2\lambda^2 n}}\bigl(1+ C(\theta)\d\bigr).
\]
We may also write for any {$d \geq0$}
\[
\int_{\afrac{d}{\lambda\sqrt n}}^\infty y \mathrm{e}^{\a d-\psi_n y}\mathrm{e}^{-\sfrac
{y^2}{2}}\,\mathrm{d}y
\leq\int_{\afrac{d}{\lambda\sqrt n}}^\infty y \mathrm{e}^{-\sfrac
{y^2}2}\,\mathrm{d}y\leq
\mathrm{e}^{-\afrac{d^2}{2\lambda^2 n}}.
\]
Hence,
\[
\sqrt{2\pi} \mathrm{e}^{\a d} \bigl|I(n)\bigr| \leq2\mathrm{e}^{-\afrac{d^2}{2\lambda^2 n}}.
\]
Now we compute the second term $\mathit{II}(n)$. Denote $g(x) =\sqrt{2\pi}
(1-x^2) \phi(x)$. Then
\begin{eqnarray*}
\sqrt{2\pi} \mathit{II}(n) &=& \frac{m_3 \a\lambda}{6} \int_{\a d}^\infty
\mathrm{e}^{-x} \biggl[ g \biggl( \frac{x}{\psi_n} \biggr)- g \biggl(
\frac
{d}{\lambda\sqrt n} \biggr) \biggr]\,\mathrm{d}x
\\
&=& C\int_{\a d}^\infty \mathrm{e}^{-x} \int
_{\afrac{d}{\lambda\sqrt
n}}^{\sfrac{x}{\psi_n}} g'(y)\,\mathrm{d}y\,\mathrm{d}x
\\
&=& C\int_{\afrac{d}{\lambda\sqrt n}}^\infty g'(y) \int
_{\psi_n
y}^{\infty } \mathrm{e}^{-x}\,\mathrm{d}x \,\mathrm{d}y
\\
&=& C\int_{\afrac{d}{\lambda\sqrt n}}^\infty \mathrm{e}^{-\psi_n y}
g'(y) \,\mathrm{d}y.
\end{eqnarray*}
Hence,
\[
\sqrt{2\pi}\bigl|\mathit{II} (n)\bigr|\leq C \int_{\afrac{d}{\lambda\sqrt n}}^\infty
\mathrm{e}^{-\psi_n y} \,\mathrm{d}y= -\frac{C}{\psi_n} \mathrm{e}^{-\psi_n y} \biggl|_{\afrac{d}{\lambda\sqrt n}}^\infty
= \frac{C}{\lam\a\sqrt n} \mathrm{e}^{-\a d},
\]
and so
\[
\mathrm{e}^{\a d}\bigl|\mathit{II}(n)\bigr| = \mathrm{O} \biggl(\frac{1}{\sqrt n} \biggr).
\]
Finally,
\[
\sqrt{2\pi}\mathrm{e}^{\a d}J\leq \mathrm{o}(1)+2\mathrm{e}^{-\afrac{d^2}{2\lambda^2
n}}+\mathrm{O}\biggl(
\frac{1}{\sqrt{n}}\biggr),
\]
which shows \eqref{eq:p1} and
\[
\sqrt{2\pi}\mathrm{e}^{\a d}J= \mathrm{o}(1)+\mathrm{e}^{-\afrac{d^2}{2\lambda^2n}}\bigl(1+C(\theta)\d\bigr) +\mathrm{O}
\biggl(\frac{1}{\sqrt{n}}\biggr).
\]
{We may always take $\d=\d(n) =
\mathrm{o}(1)$. Hence \eqref{eq:p2} follows.}
\end{pf}

\end{appendix}

\section*{Acknowledgements}
The authors were supported in part by NCN Grant DEC-2012/05/B/ST1/00692.


\printhistory


\begin{thebibliography}{18}


\bibitem{AR}
\begin{barticle}[mr]
\bauthor{\bsnm{Addario-Berry},~\bfnm{Louigi}\binits{L.}} \AND
\bauthor{\bsnm{Reed},~\bfnm{Bruce}\binits{B.}}
(\byear{2009}).
\btitle{Minima in branching random walks}.
\bjournal{Ann. Probab.}
\bvolume{37}
\bpages{1044--1079}.
\bid{doi={10.1214/08-AOP428}, issn={0091-1798}, mr={2537549}}
\end{barticle}
\bptok{imsref}%
\endbibitem

\bibitem{ABM}
\begin{barticle}[mr]
\bauthor{\bsnm{Alsmeyer},~\bfnm{Gerold}\binits{G.}},
\bauthor{\bsnm{Biggins},~\bfnm{J.~D.}\binits{J.D.}} \AND
\bauthor{\bsnm{Meiners},~\bfnm{Matthias}\binits{M.}}
(\byear{2012}).
\btitle{The functional equation of the smoothing transform}.
\bjournal{Ann. Probab.}
\bvolume{40}
\bpages{2069--2105}.
\bid{doi={10.1214/11-AOP670}, issn={0091-1798}, mr={3025711}}
\end{barticle}
\bptok{imsref}%
\endbibitem

\bibitem{ADM}
\begin{bmisc}[auto:STB|2014/01/06|10:16:28]
\bauthor{\bsnm{Alsmeyer},~\bfnm{G.}\binits{G.}},
\bauthor{\bsnm{Damek},~\bfnm{E.}\binits{E.}} \AND
\bauthor{\bsnm{Mentemeier},~\bfnm{S.}\binits{S.}}
\bhowpublished{Precise tail index of fixed points of the two-sided smoothing transform. \arxivurl{arXiv:1206.3970v1}}.
\end{bmisc}
\bptok{imsref}%
\endbibitem

\bibitem{AM2}
\begin{barticle}[mr]
\bauthor{\bsnm{Alsmeyer},~\bfnm{Gerold}\binits{G.}} \AND
\bauthor{\bsnm{Meiners},~\bfnm{Matthias}\binits{M.}}
(\byear{2012}).
\btitle{Fixed points of inhomogeneous smoothing transforms}.
\bjournal{J. Difference Equ. Appl.}
\bvolume{18}
\bpages{1287--1304}.
\bid{doi={10.1080/10236198.2011.589514}, issn={1023-6198}, mr={2956046}}
\end{barticle}
\bptok{imsref}%
\endbibitem

\bibitem{AM}
\begin{barticle}[mr]
\bauthor{\bsnm{Alsmeyer},~\bfnm{Gerold}\binits{G.}} \AND
\bauthor{\bsnm{Meiners},~\bfnm{Matthias}\binits{M.}}
(\byear{2013}).
\btitle{Fixed points of the smoothing transform: Two-sided solutions}.
\bjournal{Probab. Theory Related Fields}
\bvolume{155}
\bpages{165--199}.
\bid{doi={10.1007/s00440-011-0395-y}, issn={0178-8051}, mr={3010396}}
\end{barticle}
\bptok{imsref}%
\endbibitem

\bibitem{bassetti}
\begin{barticle}[mr]
\bauthor{\bsnm{Bassetti},~\bfnm{Federico}\binits{F.}} \AND
\bauthor{\bsnm{Ladelli},~\bfnm{Lucia}\binits{L.}}
(\byear{2012}).
\btitle{Self-similar solutions in one-dimensional kinetic models: A probabilistic view}.
\bjournal{Ann. Appl. Probab.}
\bvolume{22}
\bpages{1928--1961}.
\bid{doi={10.1214/11-AAP818}, issn={1050-5164}, mr={3025685}}
\end{barticle}
\bptok{imsref}%
\endbibitem

\bibitem{BDMM}
\begin{barticle}[mr]
\bauthor{\bsnm{Buraczewski},~\bfnm{Dariusz}\binits{D.}},
\bauthor{\bsnm{Damek},~\bfnm{Ewa}\binits{E.}},
\bauthor{\bsnm{Mentemeier},~\bfnm{Sebastian}\binits{S.}} \AND
\bauthor{\bsnm{Mirek},~\bfnm{Mariusz}\binits{M.}}
(\byear{2013}).
\btitle{Heavy tailed solutions of multivariate smoothing transforms}.
\bjournal{Stochastic Process. Appl.}
\bvolume{123}
\bpages{1947--1986}.
\bid{doi={10.1016/j.spa.2013.02.003}, issn={0304-4149}, mr={3038495}}
\end{barticle}
\bptok{imsref}%
\endbibitem

\bibitem{DZ}
\begin{bbook}[mr]
\bauthor{\bsnm{Dembo},~\bfnm{Amir}\binits{A.}} \AND
\bauthor{\bsnm{Zeitouni},~\bfnm{Ofer}\binits{O.}}
(\byear{1998}).
\btitle{Large Deviations Techniques and Applications},
\bedition{2nd} ed.
\bseries{Applications of Mathematics (New York)}
\bvolume{38}.
\blocation{New York}:
\bpublisher{Springer}.
\bid{mr={1619036}}
\end{bbook}
\bptok{imsref}%
\endbibitem

\bibitem{DL}
\begin{barticle}[mr]
\bauthor{\bsnm{Durrett},~\bfnm{Richard}\binits{R.}} \AND
\bauthor{\bsnm{Liggett},~\bfnm{Thomas~M.}\binits{T.M.}}
(\byear{1983}).
\btitle{Fixed points of the smoothing transformation}.
\bjournal{Z. Wahrsch. Verw. Gebiete}
\bvolume{64}
\bpages{275--301}.
\bid{doi={10.1007/BF00532962}, issn={0044-3719}, mr={0716487}}
\end{barticle}
\bptok{imsref}%
\endbibitem

\bibitem{F}
\begin{bbook}[mr]
\bauthor{\bsnm{Feller},~\bfnm{William}\binits{W.}}
(\byear{1971}).
\btitle{An Introduction to Probability Theory and Its Applications. {V}ol. {II}.}
\bseries{Second Edition}.
\blocation{New York}:
\bpublisher{Wiley}.
\bid{mr={0270403}}
\end{bbook}
\bptok{imsref}%
\endbibitem

\bibitem{Go}
\begin{barticle}[mr]
\bauthor{\bsnm{Goldie},~\bfnm{Charles~M.}\binits{C.M.}}
(\byear{1991}).
\btitle{Implicit renewal theory and tails of solutions of random equations}.
\bjournal{Ann. Appl. Probab.}
\bvolume{1}
\bpages{126--166}.
\bid{issn={1050-5164}, mr={1097468}}
\end{barticle}
\bptok{imsref}%
\endbibitem

\bibitem{HS}
\begin{barticle}[mr]
\bauthor{\bsnm{Hu},~\bfnm{Yueyun}\binits{Y.}} \AND
\bauthor{\bsnm{Shi},~\bfnm{Zhan}\binits{Z.}}
(\byear{2009}).
\btitle{Minimal position and critical martingale convergence in branching random walks, and directed polymers on disordered trees}.
\bjournal{Ann. Probab.}
\bvolume{37}
\bpages{742--789}.
\bid{doi={10.1214/08-AOP419}, issn={0091-1798}, mr={2510023}}
\end{barticle}
\bptok{imsref}%
\endbibitem

\bibitem{JC1}
\begin{barticle}[mr]
\bauthor{\bsnm{Jelenkovi{\'c}},~\bfnm{Predrag~R.}\binits{P.R.}} \AND
\bauthor{\bsnm{Olvera-Cravioto},~\bfnm{Mariana}\binits{M.}}
(\byear{2010}).
\btitle{Information ranking and power laws on trees}.
\bjournal{Adv. in Appl. Probab.}
\bvolume{42}
\bpages{1057--1093}.
\bid{doi={10.1239/aap/1293113151}, issn={0001-8678}, mr={2796677}}
\end{barticle}
\bptok{imsref}%
\endbibitem

\bibitem{JC2}
\begin{barticle}[mr]
\bauthor{\bsnm{Jelenkovi{\'c}},~\bfnm{Predrag~R.}\binits{P.R.}} \AND
\bauthor{\bsnm{Olvera-Cravioto},~\bfnm{Mariana}\binits{M.}}
(\byear{2012}).
\btitle{Implicit renewal theory and power tails on trees}.
\bjournal{Adv. in Appl. Probab.}
\bvolume{44}
\bpages{528--561}.
\bid{doi={10.1239/aap/1339878723}, issn={0001-8678}, mr={2977407}}
\end{barticle}
\bptok{imsref}%
\endbibitem

\bibitem{JC3}
\begin{barticle}[mr]
\bauthor{\bsnm{Jelenkovi{\'c}},~\bfnm{Predrag~R.}\binits{P.R.}} \AND
\bauthor{\bsnm{Olvera-Cravioto},~\bfnm{Mariana}\binits{M.}}
(\byear{2012}).
\btitle{Implicit renewal theorem for trees with general weights}.
\bjournal{Stochastic Process. Appl.}
\bvolume{122}
\bpages{3209--3238}.
\bid{doi={10.1016/j.spa.2012.05.004}, issn={0304-4149}, mr={2946440}}
\end{barticle}
\bptok{imsref}%
\endbibitem

\bibitem{NR}
\begin{barticle}[mr]
\bauthor{\bsnm{Neininger},~\bfnm{Ralph}\binits{R.}} \AND
\bauthor{\bsnm{R{\"u}schendorf},~\bfnm{Ludger}\binits{L.}}
(\byear{2004}).
\btitle{A general limit theorem for recursive algorithms and combinatorial structures}.
\bjournal{Ann. Appl. Probab.}
\bvolume{14}
\bpages{378--418}.
\bid{doi={10.1214/aoap/1075828056}, issn={1050-5164}, mr={2023025}}
\end{barticle}
\bptok{imsref}%
\endbibitem

\bibitem{R}
\begin{barticle}[mr]
\bauthor{\bsnm{R{\"o}sler},~\bfnm{U.}\binits{U.}}
(\byear{2001}).
\btitle{On the analysis of stochastic divide and conquer algorithms}.
\bjournal{Algorithmica}
\bvolume{29}
\bpages{238--261}.
\bnote{Average-case analysis of algorithms (Princeton, NJ, 1998)}.
\bid{doi={10.1007/BF02679621}, issn={0178-4617}, mr={1887306}}
\end{barticle}
\bptok{imsref}%
\endbibitem

\bibitem{VL}
\begin{barticle}[mr]
\bauthor{\bsnm{Volkovich},~\bfnm{Yana}\binits{Y.}} \AND
\bauthor{\bsnm{Litvak},~\bfnm{Nelly}\binits{N.}}
(\byear{2010}).
\btitle{Asymptotic analysis for personalized web search}.
\bjournal{Adv. in Appl. Probab.}
\bvolume{42}
\bpages{577--604}.
\bid{doi={10.1239/aap/1275055243}, issn={0001-8678}, mr={2675117}}
\end{barticle}
\bptok{imsref}%
\endbibitem

\end{thebibliography}
\end{document}